\documentclass[12pt,reqno]{amsart}

\usepackage{amssymb}
\usepackage{hyperref}
\usepackage[all]{xypic}

\newcommand{\F}{\mathbb{F}}

\newcommand{\Fp}{\F_p}

\newcommand{\Gal}{\text{\rm Gal}}

\newcommand{\Z}{\mathbb{Z}}
\newcommand{\ch}[1]{\textrm{char(}#1\textrm{)}}
\newcommand{\comment}[1]{}
\newcommand{\inv}{^{-1}}
\begin{document}

\title{$p$-groups have unbounded realization multiplicity}

\author{Jen Berg}
\address{Department of Mathematics, 1 University Station C1200, Austin, TX 78712-0257}
\email{jberg@math.utexas.edu}

\author[Andrew Schultz]{Andrew Schultz}
\address{Department of Mathematics, Wellesley College, 106 Central Street, Wellesley, MA 02482}
\email{andrew.c.schultz@gmail.com}

\begin{abstract}In this paper we interpret the solutions to a particular Galois embedding problem over an extension $K/F$ satisfying $\Gal(K/F) \simeq \Z/p^n\Z$ in terms of certain Galois submodules within the parameterizing space of elementary $p$-abelian extensions of $K$; here $p$ is a prime.  Combined with some basic facts about the module structure of this parameterizing space, this allows us to exhibit a class of $p$-groups whose realization multiplicity is unbounded.
\end{abstract}

\date{\today}

\maketitle

\newtheorem*{theorem*}{Theorem}
\newtheorem*{lemma*}{Lemma}
\newtheorem{proposition}{Proposition}[section]
\newtheorem{theorem}[proposition]{Theorem}
\newtheorem{corollary}[proposition]{Corollary}
\newtheorem{lemma}[proposition]{Lemma}

\theoremstyle{definition}
\newtheorem*{definition*}{Definition}
\newtheorem*{remark*}{Remark}

\parskip=10pt plus 2pt minus 2pt

\section{Introduction}

In the study of Galois extensions of a given field $F$, one is interested in determining those field extensions of $F$ whose Galois group is some given group $G$.  If $K/F$ is an extension with $\Gal(K/F) \simeq G$, we say that $K$ is a $G$-extension of $F$.  We write $\nu(G,F)$ for the number of distinct $G$-extensions of $F$ within a fixed algebraic closure of $F$.  We then write $\mathfrak{F}(G)$ for the set of fields $F$ such that $\nu(G,F) \geq 1$.  

There has been a great deal of interest in studying the relationships amongst the collections $\mathfrak{F}(G)$ for various $G$, particularly because this provides a means for studying how absolute Galois groups differ from ``random'' profinite groups.  For example, for groups $G$ and $H$ and $F$ a field, one can ask whether $F \in \mathfrak{F}(G)$ implies $F \in \mathfrak{F}(H)$; if so, then $G$ is said to automatically realize $H$.  The trivial automatic realizations are those where $H$ is a quotient of $G$.  The first non-trivial automatic realization was given by Whaples in \cite{Wh}, but the hunt for non-trivial automatic realization results began in earnest with \cite{J1}.  Since then, there have been a number papers which give non-trivial automatic realization results: \cite{B,GSS,J2,J3,M1,M2}.

In this paper, we consider a close cousin of automatic realization results.  The realization multiplicity of $G$, written $\nu(G)$, is defined as $$\nu(G) = \min_{F \in \mathfrak{F}(G)}\{\nu(G,F)\}.$$  When realization multiplicities were first studied, it was observed that a number of groups have realization multiplicity $1$.  For instance, any cyclic group has realization multiplicity $1$ since $\Gal(\bar \F_p ^{\text{sep}}/\F_p) = \hat \Z$, and it can be shown that all groups of order at most $15$ have realization multiplicity $1$.  In \cite{J2}, Jensen says that M.~Jarden asked whether all groups have realization multiplicity $1$.  Jensen goes on to answer this question in the negative by showing that $$\nu(D_4 \times \Z/2\Z) = 3 \quad \mbox{ and } \quad \nu(SD_8) = 2.$$  Here we write $D_{2^{n-1}}$ for the dihedral group of order $2^n$ ($n \geq 3$), and $SD_{2^{n-1}}$ for the quasi-dihedral group of order $2^n$ ($n \geq 4$).  Both groups are generated by elements $\sigma$ and $\tau$ subject to the relation $\sigma^{2^{n-1}} = \tau^2 = 1$, though the former group has $\tau \sigma\tau = \sigma^{-1}$ and the latter $\tau \sigma \tau = \sigma^{2^{n-2}}-1$.

Since this initial paper, much of the work in this area has focused on computing realization multiplicities for $2$-groups.  In particular, \cite[Prop.~2.3.1 and Prop.~2.3.5]{JP2} show that the only non-abelian group of order $2^n$ and exponent $2^{n-1}$ whose realization multiplicity is greater than $1$ is $QD_{2^{n-1}}$, with $\nu(SD_{2^{n-1}}) = 2.$  The question of whether or not $2$-groups have unbounded realization multiplicity was settled in this same paper, where it was shown that $\nu(D_4^n) \geq n$.  Jensen observes in \cite{J2} that ``[i]t is a little harder to obtain indecomposable $2$-groups of realization multiplicity $>$1," and indeed it seems that the aforementioned $SD_{2^{n-1}}$ are the only known examples.

There has been relatively little progress towards understanding realization multiplicities for $p$-groups when $p>2$; aside from some abelian $p$-groups with realization multiplicity $1$, the known results are
\begin{itemize}
\item $\nu\left(H(p^3)\times (\Z/p\Z)^n\right) = 1;$ 
\item $\nu\left(M(p^3)\right) = 2;$ and
\item $\nu\left(M(p^3)\times \Z/p\Z \right) = p^2-1$,
\end{itemize}
where $H(p^3)$ is the non-abelian group of order $p^3$ and exponent $p$, while $M(p^3)$ is the non-abelian group of order $p^3$ and exponent $p^2$.

The purpose of this paper is to extend the class of $p$-groups whose realization multiplicity is known to be greater than $1$, and particularly to show that $p$-groups have unbounded realization multiplicity.  The central theorem we prove is

\begin{theorem}\label{th:realization_multiplicity}
Suppose that $p$ is prime and $n$ is a positive integer, with $n \geq 2$ when $p=2$.  Let $k$ be given.  Then 
$$\nu\left(\F_p[\Z/p^n\Z]^k\rtimes \Z/p^n\Z\right) \geq p^k.$$
\end{theorem}

In fact, our methodology is to count the appearance of this group over $F$ by studying the extensions $L/F$ over a given extension $K/F$ with $\Gal(K/F) \simeq \Z/p^n\Z := G$ so that the natural short exact sequence from Galois theory is
$$\xymatrix{\F_p[\Z/p^n\Z]^k \rtimes \Z/p^n\Z \ar[r] & \Z/p^n\Z \ar[r] & 1}.$$  In other words, our theorem is really a lower bound on the number of solutions to this particular Galois embedding problem given the existence of at least one solution.  

This paper also serves as an advertisement for a new approach to studying embedding problems with factor group $\Z/p^n\Z$ and elementary $p$-abelian kernel. Often such embedding problems are studied through the lens of the Brauer group and its connection to known obstructions of embedding problems. On the contrary, we reach our results by interpreting the appearance of certain Galois groups over a field $F$ in terms of the Galois-module structure of well-known parameterizing spaces of elementary $p$-abelian extensions.  This methodology has already been used to compute automatic realizations for certain groups (see \cite{MS2,MSS2}), and the second author of this paper extends the ideas developed in this paper to give module-theoretic interpretations for solutions to any embedding problems $\hat G \to G \to 1$ for any group $\hat G$ which is an extension of $G=\Z/p^n\Z$ by an elementary $p$-abelian group in \cite{S}. In particular this allows for a precise count of solutions to these embedding problems, as well as giving new automatic realization and realization multiplicity results.

The tools used in studying these problems for fields with $\ch{K} \neq p$ were already introduced in \cite{MSS1}.  Though we could use Witt's theorem (or its progeny) on the realizability of $p$-groups over fields of characteristic $p$, we instead provide a module-theoretic interpretation in this case that echoes the case $\ch{K} \neq p$ (though, of course, Witt's theorem will still be useful for us).  We believe the small cost of translating this Galois module methodology into the characteristic $p$ setting is justified by the uniform methodology it presents us in searching for extensions with a prescribed group, particularly when it comes to parameterizing such extensions. 

The paper is structured as follows.  In section \ref{sec:known} we review established notation and known results for fields whose characteristic is different from $p$.  In section \ref{sec:artin} we carry these notions over into the characteristic $p$ setting by way of Artin-Schreier theory.  In section \ref{sec:realization} we translate common features from all these settings into the desired realization multiplicity result.

\section{Notation and Known Results}\label{sec:known}

Throughout the paper we will use $G$ to denote the group $\Z/p^n\Z$, and we let $\sigma$ be a fixed generator, where $p$ is prime and $n$ is a positive integer with $n \geq 2$ when $p=2$.  The extension $K/F$ will always be a $G$-extension.

Kummer theory tells us that if a field $K$ contains a primitive $p$th root of unity $\xi_p$, then extensions $L/K$ with $\Gal(L/K) \simeq \left(\Z/p\Z\right)^k$ are parametrized by $k$-dimensional subspaces of $J=K^\times/K^{\times p}$.  When $K$ is itself an extension of $F$ so that $\Gal(K/F) \simeq G$, then Waterhouse studied in \cite{Wa} how the $G$-action on $J$ could be used to detect those elementary $p$-abelian extensions of $K$ which were Galois over $F$.  Indeed, in the event that one chooses a cyclic $\F_p[G]$-submodule $N \subseteq J$ and $L/K$ is the corresponding extension, Waterhouse was able to compute $\Gal(L/F)$.  

He also investigates analogous questions when $\ch{K} \neq p$ but $\xi_p \not\in K$.  He shows that there is a subspace $J \subseteq K(\xi_p)^\times/K(\xi_p)^{\times p}$ which parametrizes elementary $p$-abelian extensions of $K$, and whose $\F_p[G]$-submodules correspond to elementary $p$-abelian extensions of $K$ which are additionally Galois over $F$ (the $G$-action on $J$ comes from the natural extension of the $G$-action to $K(\xi_p)$). 

In this paper we will focus on those elementary $p$-abelian extensions $L/K$ so that $\Gal(L/F) \simeq \F_p[G] \rtimes G$.  Waterhouse provides us with the following characterization of such extensions.



\begin{proposition}\label{prop:module.to.galois.translation.kummer}
Suppose that $K$ is a field with $\ch{K} \neq p$.  Let $N\subseteq J$ be an $\F_p[G]$-submodule, and let $L/K$ be the corresponding elementary $p$-abelian extension.  Then $\Gal(L/F) \simeq \F_p[G] \rtimes G$ if $N \simeq \F_p[G]$.
\end{proposition}

\section{Generalizing to characteristic $p$}\label{sec:artin}

When one investigates $p$-groups as Galois groups over fields of characteristic $p$, the indispensable tool is Witt's theorem from \cite{Wi}.  Witt's methodology was adapted to the context of embedding problems in \cite[App.~A]{JLY}, where it is shown that any central, non-split embedding problem
$$1\to \Z/p\Z \to G \to H \to 1$$ has a solution over any field of characteristic $p$.  Though this result will be useful to us in the next section, we will choose to study the appearance of certain $p$-groups over fields of characteristic $p$ in a way analogous to Waterhouse's approach for Kummer theory.  This allows us to make module-theoretic arguments which work uniformly for all fields, regardless of characteristic or the presence of roots of unity.

We remind the reader that in characteristic $p$, Kummer theory is replaced by Artin-Schreier theory, which tells us that elementary $p$-abelian extensions of a field $K$ with $\ch{K} = p$ are parametrized by $\F_p$-spaces of $$K/\wp(K), \mbox{ where }\wp(K) = \{k^p-k: k \in K\}.$$  (The action of $c \in \F_p$ on an element $\gamma \in K/\wp(K)$ is given by $\gamma \mapsto \gamma+c$.) Since all the structures in this setting are additive, we will translate all our notation over into additive form in this section.  

The Artin-Schreier parametrization is as follows. For an element $\gamma \in K/\wp(K)$, write $\rho(\gamma)$ for a root of the polynomial $x^p-x-\gamma$.  Then if $N \subseteq K/\wp(K)$, the field $L$ to which $N$ corresponds is given by $L = K(\rho(\gamma): \gamma \in N).$  There is a pairing $\Gal(L/K) \times N \to \F_p$ defined by $$\langle \tau,\gamma \rangle \mapsto \tau(\rho(\gamma))-\rho(\gamma) = (\tau-1)\rho(\gamma)$$ which is independent of the choice of root $\rho(\gamma)$, and Artin-Schreier theory tells us that this pairing is perfect; i.e., $N$ and $\Gal(L/K)$ are dual as $\F_p$-spaces.

As in the previous section, we are interested in the case where $K$ is an extension of a field $F$ with $G=\Gal(K/F) \simeq \Z/p^n\Z$, and particularly information about the Galois structure of $L/F$ when $L/K$ is an elementary $p$-abelian extension.  For instance, it is an easy exercise to show that such an extension $L/F$ is Galois over $F$ if and only if the corresponding $\F_p$-space $N \subseteq K/\wp(K)$ is actually an $\F_p[G]$-submodule.  One also has that the duality between $\Gal(L/K)$ and $N$ is $G$-equivariant, as shown in the following

\begin{lemma}
If $\ch{K}=p$ and $N\subseteq K/\wp(K)$ is the $\F_p[G]$-submodule corresponding to an elementary $p$-abelian extension $L/K$, then $N$ and $\Gal(L/K)$ are $\F_p[G]$-dual under the Artin-Schreier pairing.
\end{lemma}

\begin{proof}
Let $\tau \in \Gal(L/K)$, $\gamma \in N$, and let $\tilde \sigma$ denote a lift of the generator of $G$ to $\Gal(L/K)$. We check that $$ \langle \sigma \cdot \tau, \sigma \gamma \rangle = \sigma \langle \tau, \gamma \rangle $$ where $\sigma \cdot \tau = \sigma \tau \sigma \inv$.

We note that if $\rho(\gamma)$ is a root of $x^p - x - \gamma$ then $\tilde \sigma(\rho(\gamma))$ is a root of $x^p - x - \sigma(\gamma)$. Since the pairing is independent of the choice of $\rho(\sigma(\gamma))$, we have that

\begin{eqnarray*} \langle \sigma \cdot \tau, \sigma (\gamma) \rangle & = & \tilde \sigma \tau \tilde\sigma \inv \left( \tilde \sigma(\rho(\gamma))\right) - \tilde \sigma(\rho(\gamma)) \\
& = & \tilde \sigma \tau \left( \rho(\gamma)\right) -\tilde \sigma (\rho(\gamma))\\
& = & \tilde \sigma \left(\tau (\rho(\gamma)) - \rho(\gamma)\right) \\
& = & \tilde \sigma\left( (\tau -1)\rho(\gamma)\right) \\
& = & \sigma \langle \tau, \gamma \rangle
\end{eqnarray*} 
\end{proof}

\begin{remark*}
Waterhouse uses this same argument to show that the analogous statements are true when $\ch{K} \neq p$.
\end{remark*}

In this paper we will focus only on those elementary $p$-abelian extensions $L/K$ which correspond to free $\F_p[G]$-modules in $K/\wp(K)$.  With this in mind, we now give the generalization of Proposition \ref{prop:module.to.galois.translation.kummer} to Artin-Schreier theory.

\begin{proposition}\label{prop:module.to.galois.translation.artin.schreier}
Suppose that $K$ is a field with $\ch{K} = p$ and $\Gal(K/F) \simeq G$.  Let $L/K$ be an elementary $p$-abelian extension which is Galois over $F$, and let $N\subseteq J$ be the corresponding $\F_p[G]$-submodule.  Then $\Gal(L/F) \simeq \F_p[G] \rtimes G$ if $N \simeq \F_p[G]$.
\end{proposition}

\begin{proof}
It is shown in \cite[p.~3]{MSS2} that $\F_p[G]$-modules are self-dual, so $\Gal(L/K) \simeq \F_p[G]$ if $N \simeq \F_p[G]$.  Since the $G$-action on $\F_p[G]$ is compatible with the $G$-action on $\Gal(L/K)$, we know that $\Gal(L/F)$ is an extension of $G$ by $\F_p[G]$.  Waterhouse shows that the only such extension is $\F_p[G]\rtimes G$ (\cite[Th.~2]{Wa}).
%
%
\end{proof}

\section{Finding large realization multiplicities}\label{sec:realization}

We've seen that for any field $K$ there is a natural parametrizing space for elementary $p$-abelian extensions, and that when there exists a field $F$ so that $\Gal(K/F) \simeq \Z/p^n\Z =G$, the $\F_p[G]$-structure on this parametrizing space can be used to describe solutions to the embedding problem
$$\F_p[G]\rtimes G \to \Gal(K/F) \simeq G \to 1.$$  In this section we give a module-theoretic interpretation for the solution to the embedding problem $$(\F_p[G])^k \rtimes G \to \Gal(K/F) \simeq G \to 1,$$ and we show that the appearance of one such extension forces the appearance of many others.

Throughout this section we will use $J$ to denote the parametrizing space of elementary $p$-abelian extensions of a field $K$, and we will always write the $\F_p[G]$-action multiplicatively.

\begin{lemma}\label{le:realizing_big_semidirect_products}
Let $L/K$ be an elementary $p$-abelian extension which is Galois over $F$, and let $N \subseteq J$ be the corresponding $\F_p[G]$-submodule.  Then $\Gal(L/F) \simeq \left(\F_p[G]\right)^k \rtimes G$ if $N \simeq \left(\F_p[G]\right)^k$. 
\end{lemma}
\begin{proof}
As in the proof of Proposition \ref{prop:module.to.galois.translation.artin.schreier}, the duality of $\Gal(L/K)$ with $N$, together with the self-dual nature of $\F_p[G]$-modules, means that $\Gal(L/K) \simeq \left(\F_p[G]\right)^k$ if $N \simeq \left(\F_p[G]\right)^k$.  Hence we have that $\Gal(L/F)$ is an extension of $G$ by $\left(\F_p[G]\right)^k$, and we must show that it is the semi-direct product.

So let $\tau_i\in \Gal(L/K)$ generate the $i$th summand of $\left( \F_p[G]\right)^k$.  Since $\langle \tau_i \rangle \simeq \F_p[G]$, it corresponds to an extension $L_i/F$ with $\Gal(L_i/F) \simeq \F_p[G]\rtimes G$.  The various $L_i$ have pairwise intersection $K$ since their corresponding modules in $J$ have trivial intersection, and their compositum is $L$ since the $\tau_i$ generate $N$.  Hence $\Gal(L/F)$ is the subgroup of elements $$\left\{\prod_{i=1}^k (n_i,\sigma^{j_i}) ~:~ \mbox{res}_K(n_1,\sigma^{j_1}) = \cdots = \mbox{res}_K(n_k,\sigma^{j_k})\right\}\subseteq \prod_{i=1}^k \Gal(L_i/F).$$ But because $(n_i,\sigma^{j_i})$ restricts to $\sigma^{j_i} \in \Gal(K/F)$, this means that 
\begin{equation*}
\begin{split}
\Gal(L/F) &= \left\{\prod_{i=1}^k(n_i,\sigma^j)~:~n_i \in \langle \tau_i \rangle, 0 \leq j \leq p^n-1\right\}
\end{split}
\end{equation*} 
which is clearly isomorphic to $\left(\F_p[G]\right)^k \rtimes G$ since each $\Gal(L_i/F) \simeq \F_p[G]\rtimes G$.
\end{proof}

Our goal will be to show that the appearance of one submodule $\left(\F_p[G]\right)^k$ within $J$ forces the appearance of many other such submodules. If we can show that $J$ contains some element outside this free submodule, this result is just a matter of constructing such submodules.

\begin{proposition}\label{prop:lots_of_free}
Suppose that $W$ is an $\F_p[G]$-module which properly contains a submodule $V \simeq \left(\F_p[G]\right)^k$, and further suppose that there exists $\delta \in W \cap \ker(\sigma-1)^{p^n-1} \setminus V$.  Then $W$ contains at least $p^k$-many submodules isomorphic to $\left(\F_p[G]\right)^k$.
\end{proposition}
\begin{proof}
Let $v_i$ generate the $i$th summand of $V$ and choose an element $\vec{c} \in \F_p^k$. Define $V_{\vec{c}} = \langle \delta^{c_1}v_1,\cdots,\delta^{c_k}v_k\rangle$. We will prove that each $V_{\vec{c}}$ is isomorphic to $\Fp[G]^k$. Given that this is the case, it is not difficult to show that each $V_{\vec{c}}$ is distinct. If this were not so, then there would exist some $\vec{b} \ne \vec{c}$ such that $V_{\vec{c}} = V_{\vec{b}}$. Let $j$ be the first index where $\vec{b}$ and $\vec{c}$ differ. Then, $v_j \delta ^{c_j} v_j \inv \delta^{-b_j} = \delta^{c_i - b_i} \in V_{\vec{c}}$. Since $c_i \ne b_i$, $c_i - b_i \in \Fp^\times$ so that in fact, $\delta \in V_{\vec{c}}$. This implies that $V \subset V_{\vec{c}}$, and in fact equality holds since both are isomorphic to $\Fp[G]^k$. But then $\delta \in V$, which is impossible by the choice of $\delta$. Hence each $V_{\vec{c}}$ is distinct. Thus, it remains only to show that $V_{\vec{c}} \simeq \Fp[G]^k$.

We will show that there are no non-trivial relations amongst $\{\delta^{c_i}v_i\}$, and hence the module they generate is isomorphic to $\F_p[G]^k$.  Suppose that $\{g_i(\sigma)\} \in \F_p[G]$ satisfy $$ 1 = \prod_i (v_i \delta^{c_i})^{g_i(\sigma)} = \prod_i v_i^{g_i(\sigma)} \cdot \delta^{\sum_i c_i g_i(\sigma)}.$$ Write $g_i(\sigma) = (\sigma-1)^{r_i} \tilde{g_i}(\sigma)$, where $(\sigma - 1) \nmid \tilde{g_i}(\sigma)$, $r_i \ge 0$.  If there is some $g_i \neq 0$, then $r = \min_i \{r_i\} < p^n$.  Applying $(\sigma-1)^{p^n-1-r}$ to the given relation annihilates any term for which $r_i>r$, leaving
$$ 1 = \prod_{r_i = r} v_i^{(\sigma-1)^{p^n-1} \tilde{g_i}(\sigma)} \cdot \delta^{(\sigma - 1)^{p^n-1} (\sum_{r_i=r} c_i \tilde{g_i}(\sigma))}.$$ Since $\delta \in \ker (\sigma-1)^{p^n-1}$ we have that in fact,
$$ 1 = \prod_{r_i = r} v_i^{(\sigma-1)^{p^n - 1} \tilde{g_i}(\sigma)}.$$ This then gives a non-trivial relation amongst the $\{v_i\}$, contradicting the assumption that $V \simeq \F_p[G]^k$.
\end{proof}

We are now ready to prove Theorem \ref{th:realization_multiplicity}.

\begin{proof}[Proof of Theorem \ref{th:realization_multiplicity}]
From Proposition \ref{prop:lots_of_free} and Lemma \ref{le:realizing_big_semidirect_products}, all we need to show is that if there exists a submodule $U \subseteq J$ with $U \simeq \F_p[G]^k$, then $U \neq J$ and there exists $\delta \in J$ with $\delta \in \ker(\sigma-1)^{p^n-1}$ and $\delta \not\in U$.  In the case where $\ch{K} \neq p$ and $K$ contains a primitive $p$th root of unity, this fact is a consequence of \cite[Th.~2]{MSS1}, which states that $J$ has a cyclic summand which is not isomorphic to $\F_p[G]$.  Similarly, in the case where $\ch{K} \neq p$ and $K$ does not contain a primitive $p$th root of unity, it is known that there is a cyclic summand of $J$ which is not isomorphic to $\F_p[G]$; this is found in \cite[Th.~2]{MSS2}.  

The final case, then, is when $\ch{K} = p$.  By \cite[App.~A]{JLY}, we know that the embedding problem $$\Z/p^{n+1}\Z \to \Gal(K/F) \simeq \Z/p^n\Z \to 1$$ has a solution.  Let $L/K$ be an elementary $p$-abelian extension so that $\Gal(L/F) \simeq \Z/p^{n+1}\Z$, and let $\langle \delta \rangle \subseteq J$ be the corresponding module.  We claim that $\delta$ is not contained in any submodule isomorphic to $\left(\F_p[G]\right)^k$ within $J$, after which we can appeal to Proposition \ref{prop:lots_of_free} and Lemma \ref{le:realizing_big_semidirect_products} as before.

Suppose to the contrary that it were.  The natural projections from Galois theory $$\Gal(M/F) \to \Gal(L/F) \to \Gal(K/F) \to 1$$ then correspond to group projections $$\left(\F_p[G]\right)^k \rtimes G \to \frac{\left(\F_p[G]\right)^k \rtimes G}{H \rtimes 0} \to \frac{\left(\F_p[G]\right)^k \rtimes G}{\left(\F_p[G]\right)^k \rtimes 0} \to 1,$$ where $H$ is some index-$p$ submodule of $\left(\F_p[G]\right)^k$.  Notice in particular that $H$ must contain all elements in the image of $(\sigma-1)$. Otherwise, if $(\sigma-1)f(\sigma) \notin H$ then $\langle f(\sigma), (\sigma-1)f(\sigma)\rangle_{\F_p} \cap H = \{0\}$, and so $H$ has index at least $p^2$.  We will show that the $p^n$th power of an arbitrary element of $\left(\F_p[G]\right)^k \rtimes G$ lies in $H \rtimes 0$; this will contradict the fact that $\Gal(L/F) \simeq \Z/p^{n+1}\Z$.  

We start by noting that 
\begin{equation*}
\begin{split}
(\tau_1^{f_1(\sigma)}&,\cdots,\tau_k^{f_k(\sigma)},\sigma^j)^2 = (\tau_1^{(f_1(\sigma))(1+\sigma^j)},\cdots,\tau_k^{(f_k(\sigma))(1+\sigma^j)},\sigma^{2j}).
\end{split}
\end{equation*}
By induction we have 
\begin{equation*}
\begin{split}
(\tau_1^{f_1(\sigma)}&,\cdots,\tau_k^{f_k(\sigma)},\sigma^j)^{p^n} = (\tau_1^{(f_1(\sigma))(N(\sigma))},\cdots,\tau_k^{(f_k(\sigma))(N(\sigma))},1)
\end{split}
\end{equation*} where $N(\sigma) = 1 + \sigma^j + \cdots + (\sigma^j)^{p^n-1} \equiv (\sigma^j-1)^{p^n-1} \mod{p}$.  Since this element is clearly divisible by $\sigma-1$, the given element is in $H \rtimes 0$ as claimed.  
\end{proof}

\end{document}